\documentclass[letterpaper,journal]{IEEEtran}
\usepackage[utf8]{inputenc}
\usepackage[T1]{fontenc}
\usepackage[british]{babel}
\usepackage{amsmath,amsfonts}
\usepackage{graphicx}
\usepackage{cite}
\usepackage[svgnames]{xcolor}
\usepackage{tikz}

\usepackage[hidelinks]{hyperref}

\newcommand{\orcidID}[1]{
    \unskip\hspace{0.1em}\href{https://orcid.org/#1}{\raisebox{-0.12\height}{\includegraphics[height=0.8em]{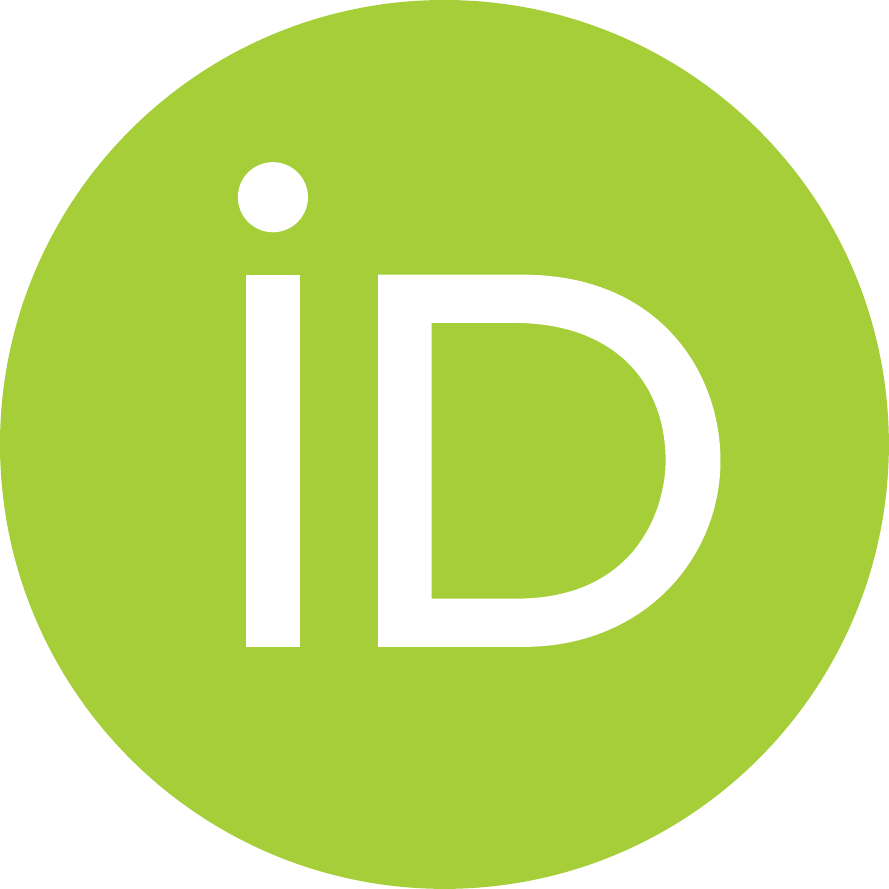}}}%
}

\begin{document}

\title{Transforming First-Year Calculus Teaching for Engineering Students -- Blocks with Field Specific Examples, Problems, and Exams}

\author{René Bødker Christensen\orcidID{0000-0002-9209-3739}, Bettina Dahl\orcidID{0000-0002-1340-3939}, and Lisbeth Fajstrup\orcidID{0000-0003-4936-1176}
    \thanks{R.B. Christensen and L. Fajstrup are with the Department of Mathematical Sciences, Aalborg University. (e-mails: \texttt{rene@math.aau.dk}, \texttt{fajstrup@math.aau.dk})}
    \thanks{B. Dahl is with the Department of Planning, Aalborg University. (e-mail: \texttt{bdahls@plan.aau.dk})}
}

% The paper headers
\markboth{Article Preprint}%
{Christensen \MakeLowercase{\textit{et al.}}: Transforming First-Year Calculus Teaching for Engineering Students}

\IEEEpubid{
    \parbox{\textwidth}{\rule{0pt}{8mm}\itshape This work has been submitted to the IEEE for possible publication. Copyright may be transferred without notice, after which this version may no longer be accessible.}
}
% Remember, if you use this you must call \IEEEpubidadjcol in the second
% column for its text to clear the IEEEpubid mark.

\maketitle

\begin{abstract}
  \emph{Contribution:} We demonstrate that it is feasible to include field specific problems in introductory mathematics courses to motivate engineering students. This is done in a way that still allows large parts of the course to be common to all students, ensuring economic viability.
  
  \emph{Background:} Many first-year engineering students perceive mathematics courses as being too abstract and far from their chosen study programme. This may lead to a lack of motivation and effort, thus decreasing the learning outcomes.
  
  \emph{Intended outcomes:} That engineering students recognize that the calculus and linear algebra courses are relevant for their future work within their specific field. This is intended to improve their learning outcomes.
  
  \emph{Application design:} The courses have been restructured into smaller subunits, each of which has a corresponding workshop treating a real-world problem from the specific field of a given group of students. These workshops are developed in consultation with the relevant fields of study, and they are given a prominent role in the course for instance by forming the basis for the oral exams.
  
  \emph{Findings:} Based on the feedback from students, we find that inclusion of field specific problems does help to highlight the applicability and importance of mathematics in engineering. When implementing such a solution, however, there are a number of challenges to keep in mind.
\end{abstract}


\section{Introduction}
\IEEEPARstart{S}{uccessful} teaching is as much about motivating students to learn as it is about providing the `right' material.

Generating this motivation can be especially challenging in introductory mathematics courses for engineering students. Various studies report that first-year engineering students are not learning mathematics the same way as mathematics students. Engineering students often complain about having to study calculus, but integrating mathematics education into the engineering field may help this student group see the use of mathematics in engineering.
Even though tools from linear algebra and calculus are vital in solving real world engineering problems, many students are discouraged by the generality and abstract nature of classical mathematical teaching. They yearn for an answer to the questions `where will I use this?' and `what's in it for me?'~\cite{Alpers13,Harterich12,Maull00,Bingolbali07,Rensaa14,Baillie00,Pepin2021}.

In this work, we evaluate a recent change in the structure of the calculus and linear algebra courses at Aalborg University (AAU) in Denmark. These courses are followed by all engineering students, primarily on their first year of study.
Both courses follow a similar structure.
To simplify the exposition, we will therefore focus on calculus in the descriptions and examples given in this work. We note, however, that the evaluation in Section~\ref{sec:evaluation} treats both courses, not only calculus.
\IEEEpubidadjcol

\subsection{Institutional setting}\label{sec:inst-sett}
The education at AAU is based on the principles of problem based learning (PBL). That is, at any given semester, the students are organized into small groups, writing a project based on PBL principles. The project work takes up half the semester with courses supporting the work.
Each group has a supervisor who facilitates the group, but the students are supposed to be self-organizing and responsible for their own learning \cite{askehave15}. 
These groups are usually self-chosen for the whole semester, but may also be formed administratively, especially on early semesters.

The project groups naturally carry over to other parts of university programmes, meaning that exercises for courses are typically solved in these groups as well.
That is, there is no need to form \emph{ad hoc} groups in teaching situations. Instead, teachers simply piggyback on the existing project groups.

Another by-product of the group-based education is the continuous feedback from students via the \emph{steering group} meetings. Each programme typically has two or three such meetings during the semester, and at these meetings, each group sends a representative to discuss any issues that may have arisen with their courses, projects, or practicalities with IT, maintenance etc.
The steering group meetings are typically organized by a \emph{semester coordinator}, and teachers are invited to participate as well.
These meetings are devoted to continual evaluation of all parts of the semester, ensuring that potential issues can be addressed -- and hopefully rectified -- quickly.

It is worth pointing out that each semester consists of a rather fixed `package' of project and courses. Students do not freely choose courses. Instead, semesters are designed in such a way that some of the courses are intended to support the projects.

The organization in groups described above means that there are certain design choices in the calculus and linear algebra courses described in this work that are natural and easy to implement at AAU, but may require more effort at other institutions.

To understand the complexity of structuring the first-year mathematics courses, one must also know that AAU has campi in three different cities: Aalborg, Copenhagen, and Esbjerg. The three cities are in opposite ends of the country with 3--5 hours of train time between them.
While Aalborg can be considered the main campus, the calculus and linear algebra courses are also followed by students in the two other cities.

\subsection{Research methodology}\label{sec:research-methodology}
The roll-out of the new course structure happened at the same time as SARS-CoV-2 (i.e. the virus causing COVID-19) spread across the world. In Denmark, this caused large parts of society to be shut down in an attempt to slow down transmission of the virus. Naturally, this also affected universities, and there were long periods of time where students and both academic and administrative staff were not allowed on campus. Courses had to be digitalized `on-the-fly', and many students were forced to spend much of their time at home. See~\cite{sstCovid} for an overview of restrictions in Denmark over the course of the pandemic.

The changes caused by COVID-19 are very likely to have had a negative impact on student well-being and student motivation, which is in turn likely to affect dropout rates and overall performance.
We have no way to control for these variables, however, so evaluating the new structure based on grades, dropout rates, and the like would not give meaningful results.

Instead, we will perform a qualitative analysis of the structure based on comments given by the students during the semester. More precisely, we will analyse the minutes of the steering group meetings described in Section~\ref{sec:inst-sett}, as these minutes also contain the students' thoughts on the calculus and linear algebra courses and their structure. This will help in gauging the students' perception of the revised structure and the use of programme-specific workshops.

In the qualitative analysis, the comments found in the minutes will be grouped based on category and sentiment -- each such combination referred to as a \emph{code}. This will take place as an inductive editing approach since we have only a few a priori codes but mainly find the final codes based on our interpretation of the meanings of these comments~\cite{robson2002real}.
Our analysis is based on all the steering group minutes we have access to from the autumn semester 2021 and the spring semester 2022. This amounts to a total of 57 meeting minutes spread across 11 study programmes.

\section{Course structure}
To allow the study programmes to tailor the course for their specific needs, they select among topics predefined by the Department of Mathematical Sciences. In this structure, each topic corresponds to a subunit, which we refer to as a \emph{block}.

\subsection{Blocks}
The initial block treats multivariate functions and their role in calculus. This includes the concept of partial derivatives and the basic principles of optimization. This is mandatory for all students following the calculus course, and to cover the required topics, this block is larger than the others, comprising six lectures. A lecture would typically consist of a four-hour session in which the first half is the actual lecture, and the second half consists of exercises where the students sit with their group members while the lecturer and potential teaching assistants walk from group to group assisting them in solving the exercises. 

The remaining seven blocks all consist of two lectures, and focus on a more specialized topic within the curriculum such as differential equations, complex numbers, or the Laplace transform.
Each study programme will only follow three of these blocks, and the choice of blocks is made by the study board responsible for the given programme. That is, from the point of view from the study board, the smaller blocks can be considered electives. For the individual student, on the other hand, the three  blocks (and thereby the curriculum) have been chosen beforehand and are mandatory.

In principle, certain blocks could have another block as a prerequisite. For instance, one might imagine a block focused on solving boundary value problems, which would then build upon a block on differential equations. For the calculus course at AAU, the content does not have such dependencies, and all the elective blocks rely only on the first, mandatory block. This gives more flexibility for the study programmes, but carries a risk of the course appearing as small disjoint pieces.
To counteract this, the examiner is given a clear, coordinating role throughout the course. We will return to the role of the examiner in Section~\ref{sec:exam}.

At AAU, the lectures in the elective blocks are livestreamed to the group rooms of the students followed by on-site exercises.
Using livestreams is not vital to the structure described in this work, however, and other educational institutions seeking to implement a similar course structure can opt for online or on-site lecturing as they see fit.

\subsection{Workshops}
Once a study programme has chosen the blocks to make up their calculus course, academic staff from the Department of Mathematical Sciences will meet with representatives from the given field of study to discuss potential problems that could form a basis for each of the four workshops.

A suitable problem for a specific block will need to satisfy two criteria. First, it should be something the students might see later in their studies or something that is currently worked on by researchers within the given field of study. Secondly, it must somehow rely on the mathematical concepts and results contained in the given block.
Note, however, that the problem itself can be much more difficult than what the students are currently capable of solving. The main purpose is to provide motivation, and to show that the calculus of the given block can be used as a tool to solve actual problems within the field.
The actual exercise problems that students will solve when doing the workshop will, for the most part, only be small `toy examples' that explore certain parts of the grand solution to the problem as a whole. Moreover, the design of the workshop ensures that  solving these exercises will require essential parts of the mathematics concepts and skills. More on choosing problems will be given in Section~\ref{sec:example-workshops}.

Finally, after developing the material for the workshop, a video introduction (\textasciitilde 20--30 min.) is made. The purpose of this video introduction is to frame the problem at hand and to present the exercise problems that the students are supposed to solve.
Thus, a workshop typically proceeds as follows: The students first meet with their group to see the video introduction and then move on to solve the exercises in the workshop. During this last part, help will be provided by teaching assistants as well as the teacher who will eventually conduct the exam.

\subsection{Exam}\label{sec:exam}
The course is evaluated through an oral exam based on the four workshops for each study programme. That is, at the exam each student draws at random one of the four workshops that the student has solved during the course. The student is then given approximately five minutes to present exercises of their choice from the given workshop. The students are encouraged to prepare and practice this presentation before the exam. Afterwards, the examiners have approximately five minutes to ask follow-up questions or other questions to the curriculum in general. Thus, including grade deliberations and other exam logistics, the examination lasts 15 minutes per student.

In Denmark, exams are usually assessed by two examiners to secure reliability and openness. One examiner is well-known to the students through being present when the students work on the workshops and through either giving some of the lectures, or being a TA in at least one of the blocks followed by these students.
The co-examiner will in most cases be a researcher from the students' subject area. The reason for this is two-fold. First, it helps to keep the evaluation in line with the expectations in a given field of study.
Second, it gives credibility to the usefulness of the course in the sense that the students see that their `own teachers' are involved in the course as well.

Note here that students do not hand in material in relation to the workshops, so the grade is based purely on the oral presentation given by the individual student.

\section{Scheduling and administration}\label{sec:sched-admin}
Splitting the course into smaller blocks that are to be followed by students from different study programmes has meant that the courses have had to be scheduled in a different way. Formerly, we would have several parallel calculus courses where the study programmes were divided between these in order to fully utilize the capacity of the lecture halls. Each course would then have one or two teachers assigned to it. This means that it would be clear to each student, who was responsible for the course and the exam.

With the increased flexibility in the curriculum, we no longer have the same flexibility in planning the parallel sessions. We now have to take into account the many different combinations of blocks, each of which is likely to have a small number of students (in some cases as low as \textasciitilde 20 students).
Running these as separate courses is not economically sustainable, and would entail many duplicate lectures. For instance the block on differential equations is followed by thirteen study programmes in the autumn, but only one of these combines this with the block on optimization, only two combine it with curves in space, and so forth.

To avoid such duplication of work, we now organize the course schedule based on the blocks, and then add study programmes as appropriate. As an example, the previously mentioned block on differential equations is conducted only once for all thirteen programmes at the same time. The students are then split up again, and regrouped based on the next block. That is, each group of students receives the chosen combination of blocks, but they will share lectures and exercise sessions with different programmes for each block.
In addition, different blocks may have different teachers assigned to them.
Note that even though the lectures may be shared by many different study programmes, they do \emph{not} share workshops since the workshop problems are tailored to each individual programme.

In Table~\ref{tab:electives} on page~\pageref{tab:electives}, we give an overview of the elective blocks of the calculus course in the autumn of 2021. This illustrates the complexity of the schedule across study programmes and across campi.

\begin{table*}[bt]
  \centering
  \includegraphics{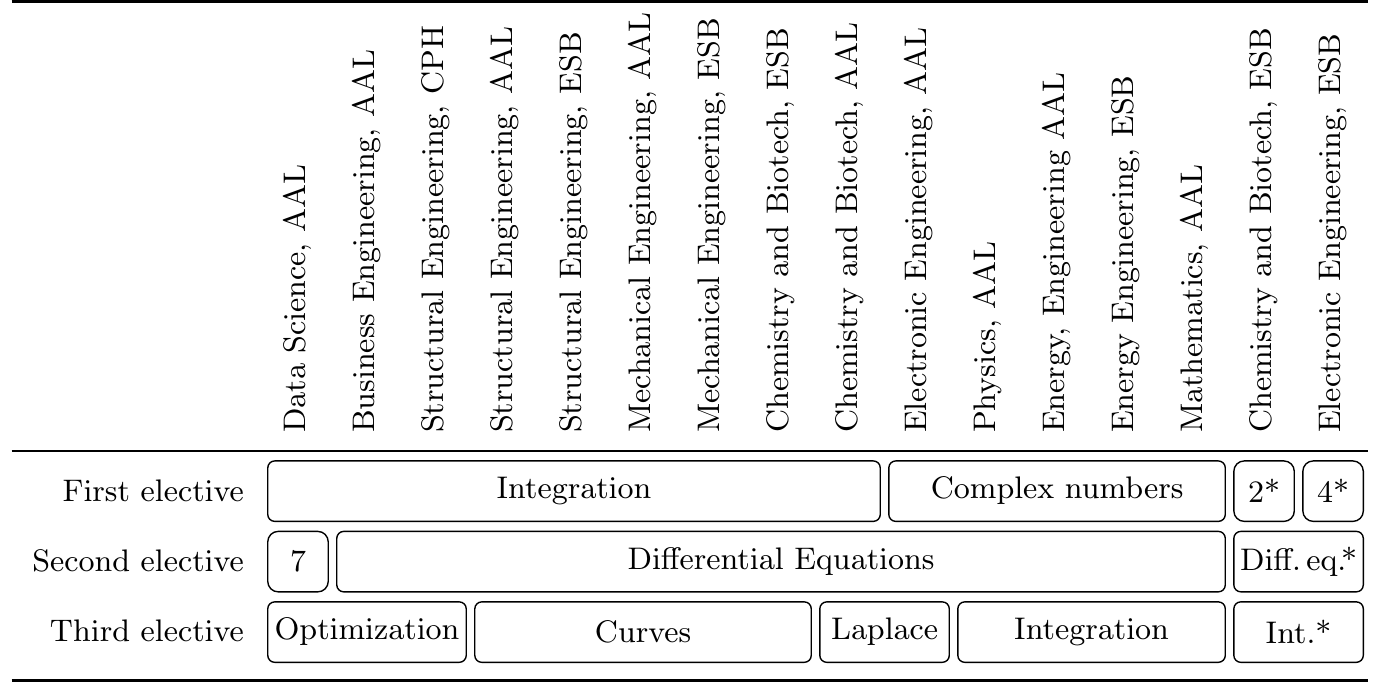}
  \caption{Autumn plan for electives in calculus for all programmes in Aalborg (AAL), Copenhagen (CPH), and Esbjerg (ESB).
      Each column corresponds to a study programme and shows the elective blocks that they follow in chronological order.
      We use general, descriptive programme titles rather than the official titles. The blocks marked 2, 4, and 7 refer to `Curves', `Complex numbers', and `Topics in Computer Science', respectively.
  Blocks marked with an asterisk (*) are conducted in English.}
  \label{tab:electives}
\end{table*}

\section{Example workshops}\label{sec:example-workshops}
In this section, we will give examples of the types of problems that have been used as motivating examples in the workshops.

\subsection{Adjustment of measurements}
All study programmes follow the first block covering partial derivatives and simple optimization, but they have different workshops.
As an example, we focus on the land surveyors, who will need to know how to handle geographic measurements subject to precision errors. In the workshop, they are given measurements of three sides in a right triangle. Their goal is to estimate the true values of these sides through least squares which for land surveyor is called adjustment. The triangle is known to be right angled and hence the estimated true values should satisfy the Pythagorean equation while the sum of the correction terms squared is minimized. Similar problems are used for other workshops for the land surveyors, e.g. Lagrange optimization and also in a block on orthogonal projection in linear algebra.

Solving the problem in this workshop entails computing a gradient and a Hessian matrix in order to determine the global minimum, leading to an analytical solution to the problem.
In the last part of the workshop, they are given a soft introduction to Newton's method via a Taylor approximation of the error function from the first part or the workshop. This is accompanied by a MATLAB-script to perform the computation as well as to compare the analytical and numerical results.

\subsection{Deflection in structural engineering}
The block on differential equations is -- among others -- followed by study programmes that can broadly be described as structural engineering.
Here, the workshop is based on the problem of determining the deflection of a bridge deck when exerted to a rhythmic force, e.g. a person walking on the bridge.
For these students, it is clear that computing such a deflection is vital for determining suitable materials in the construction of the bridge.

In the workshop, the students consider a simplified two-dimensional model of the bridge, where the deflection in question can be described by an inhomogeneous second-order ordinary differential equation.
The simplified setting does not seem to discourage students, though, since it is clear that to solve the problem in a realistic setting, one must first be able to understand and solve it in a setting with less complexity.

The workshop is then split into three sub-problems. First, the students are asked to solve the homogenous equation, finding the characteristic polynomial in the process.
Afterwards, they are guided through determining a particular solution to the inhomogeneous equation, and then combine this with the first part to obtain the complete solution.
In the final part, they are asked to solve an initial value problem, and then use MATLAB to produce a plot of the deflection over time. For this, they are given an almost finished MATLAB-script, where they only need to fill in certain values that they computed earlier during the workshop.

\subsection{Blood flow and muscle volume}
The health technology students follow the block on integration.
For their workshop, it was decided to consider two different problems related to human health.

The first is to determine the total flow of blood through a vein with a given radius and its relation to blood pressure. Using the Hagen-Poiseuille equation, the velocity of flow can be described as the solution to a differential equation, and the students are given this solution in the introduction. Determining the total flow then corresponds to integrating the velocity of flow over a cross section of the vein. The first exercise problems revolve around this as well as computing the flow in the case of a more advanced cross section.

The second problem in this workshop is the estimation of muscle volume when scanning results provide a finite number of cross sectional images along the muscle in question. The students are then guided through this estimation using truncated cones in the setting where each cross section is circular. Afterwards, they consider an example with rectangular cross sections to see the approach in more general settings.
In solving these exercise, they will use solids of revolution as well as triple integrals.

To conclude the workshop, they are asked to consider a muscle with circular cross sections given by a specific function. We provide a MATLAB-script computing the true volume, the estimated volume using truncated cones, and the estimated volume using cylinders. The students are asked to sketch the shapes whose volumes are computed and compare the results of the three calculations. Finally, they are asked to experiment with the parameters to confirm that the precisions of both estimations increase as the number of available cross sections grows.

\section{Evaluation}\label{sec:evaluation}
As described in Section~\ref{sec:research-methodology}, we will use the minutes from the steering group meetings to evaluate the student perception of the new course structure.

To aid the validity of our analysis of the minutes, we used multiple interpreters and explication of the procedures~\cite[pp.\,207--208]{Kvale96}. In terms of multiple interpreters, our procedure was as follows.
In the first step of the analysis, the first author initially analysed all the minutes, coding each utterance depending on if the utterance belonged to the category \emph{Workshop}, \emph{Teacher}, \emph{Structure}, or \emph{General} comments, respectively. The utterances were also tagged as being \emph{Positive}, \emph{Balanced}, or \emph{Negative}.
In the second step, a random programme was picked for each course using~\cite{randomOrg}, and the minutes corresponding to these programmes were then coded separately by the two other authors to check for coder reliability.
The next steps consisted of a mix of an arithmetical and a dialogical approach, aiming to reach agreement.

In the first round of samples, we discovered that the four categories were not sufficiently well-defined, and that it was often difficult to decide if an utterance was to be placed in \emph{Structure} or \emph{General}. In addition, we had left the splitting of utterances up to the individual coder, meaning that it was difficult to compare overall agreement. As an example of the latter challenge, one comment was \emph{``Block~4 was difficult to relate to from a structural engineering perspective. We refer to the block~1 workshop, which was really good in that regard.''}. One coder kept both sentences as one utterance and deemed it balanced, one split it in two being negative and positive, respectively, and one kept them together, but marked it as negative \emph{and} positive.

In order to improve coder reliability, we decided to do two things. First, reduce the categories to \emph{Workshop}, \emph{General}, and \emph{Teacher}. These were defined more precisely as
\begin{LaTeXdescription}
  \item[Workshop] Primarily concerning a workshop or its topic, video introduction etc.
  \item[General] Primarily concerning curriculum, course structure, exam etc.
  \item[Teacher] Primarily concerning the teacher, but not caused by workshops, course structure or the like.
\end{LaTeXdescription}
Second, we decided that other coders would be given the division of utterances used by the initial coder to circumvent the second type of problem seen in the first round of samples.

With these changes, a second round of samples -- different from the first round of samples -- gave an agreement of $72\%$ with respect to sentiment, $56\%$ with respect to category, and $40\%$ with respect to both at the same time. However, for every utterance at least two coders agreed on both category and sentiment.
In addition, we note that any disagreement regarding the sentiments were between balanced and negative or balanced and positive. That is, we never saw a coder perceiving an utterance as positive while another perceived it as negative.

After computing agreement percentages, we discussed the cases that we had coded differently in order to figure out why we disagreed and to obtain consensus on the most appropriate code.
During this follow-up discussion, we discovered that part of the disagreement in the three categories could be explained by hands-on experience in teaching the course. The two authors who have taught the courses tended to agree more on the categories compared to the author who only knows the structure from the general description. That is, it seems that -- unsurprisingly -- the teachers and the students share a common context that may be difficult for others to extract.

As an example of disagreement regarding sentiments, \emph{``A repetition of the expectations at the exam is wanted''} was marked as balanced by two coders and as negative by the other coder. The latter coder had put emphasis on the wish for reiterating information, suggesting that this could indicate that the initial exam information was dissatisfactory. The two other coders, on the other hand, had focused more on the wish itself being neither positive, nor negative. Namely, it is entirely possible that they simply want to make sure that they have not missed something regarding the exam.

\begin{table}[bt]
  \centering
  \includegraphics[width=\linewidth]{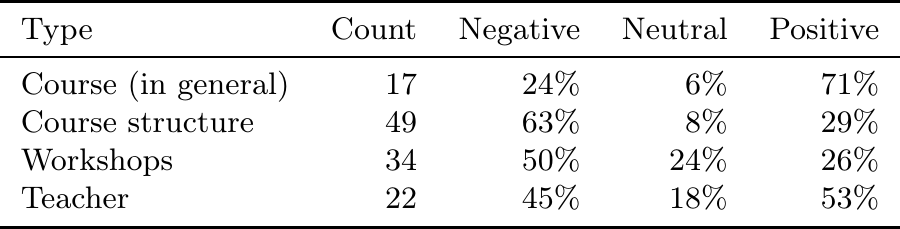}
  \caption{Distribution of student comments sorted by codes.}
  \label{tab:sentiments}
\end{table}

An overview of the comments sorted by their final codes is given in Table~\ref{tab:sentiments}.
We note here that as the focus of the steering group meetings is primarily the issues that have arisen during the semester, it is not uncommon that a majority of comments would be considered negative. If students are satisfied with a certain course, the meeting will typically move on relatively quickly, thereby devoting the time to more pressing issues.

We now proceed to analyse the overall tendencies in the comments in greater detail. Here, we are primarily interested in the comments that address the course structure and the workshops themselves.

\subsection{Workshop topics}
The health technology students comment that it is \emph{``Nice with workshops -- here you see what you need linear algebra for, quite concretely.''} They do, however, stress the importance of being able to reach and ask the responsible teacher during the workshop sessions, \emph{``otherwise one would not [gain] much from it''}.

In the initial minutes from the students in structural engineering, the motivating problems are mentioned explicitly: \emph{``During workshops, [the topics] are pinned to something tangible---works well''}. Thus, it seems like the students initially recognize the effort put into connecting terms from calculus to applications in their own field of study. At the following meeting, however, they comment that the course \emph{``\textellipsis can sometimes be difficult to relate to, and [to see] how to apply it in other contexts.''}
This could indicate that even though the students see calculus used in context during the workshops, they still find it difficult to relate to the course as a whole.

A point of critique appearing in the minutes from health technology is that their workshops had a tendency to be too extensive in comparison to the amount of time scheduled for it. In addition, they request additional resources to be allocated to have more teaching assistants available during the scheduled timeslot.
Along the same lines, the students from chemistry and biology state that \emph{``The topic worked on during the workshop seemed to be outside the curriculum[\textellipsis]''}.
This highlights a challenge regarding the choice of motivating problem. On the one hand, the problem must be as close as possible to the field of the students, but at the same time, it should not act as a hindrance to solve the underlying calculus problems.
Striking this balance can be difficult, and our experience is that the first workshops created under this new structure will need a couple of revisions before reaching the appropriate extent and level of difficulty. After this initial tuning, the subsequent workshops will be easier to develop.

\subsection{Course structure}
Some of the negative comments regarding the course structure were anticipatory. That is, rather than commenting on what \emph{has happened} in the course, they comment on their expectation of what \emph{will happen}. These comments all relate to the streaming of the elective blocks, stating that students are concerned with the prospect of online teaching.
We see this at the first meeting of 2022 for the Danish programmes in Esbjerg, where they say \emph{``We are happy that it is not yet virtual''}. Interestingly, at their last meeting where the streamed lectures have taken place, their sentiment has changed: \emph{``[The] online part went well compared to last semester.''}

A benefit of the streamed lectures is that students can watch it (or parts of it) again, which we also see in the minutes. For instance, we find comments such as
\emph{``The students like that videos stay online from the digital lectures, so that one has a resource to look back on''} and
\emph{``Lectures have been recorded, which has been nice''}.

A reoccurring theme in many of the minutes is the number of teachers involved in the course. We found comments such as
\emph{``[We] do not understand why we need different lecturers''}, 
\emph{``Is it possible to reduce the number of lecturers[\textellipsis]. Differences in presentation of material and ability to help varies a lot''}, and 
\emph{``Many different teachers, even within individual topics, which the students sometimes find confusing''}.
Part of the problem seemed to stem from us not providing enough (or sufficiently clear) information about the course structure and the assigned teachers at the beginning of the course. As a result, we now spend more time during the first lecture to explain the structure of the course, how it is organized in our learning management system, and which teacher will be responsible for which parts.
The latter is especially important with respect to the examiner; it should be clear to students who will conduct their exam. See Section~\ref{sec:exam} about the examiners role during the course.

\section{Conclusion}
In this work, we have described the restructuring of the calculus and linear algebra courses at Aalborg University in Denmark. Aiming to motivate students and to give them examples closer to their own field of study, the restructuring divides the courses into blocks, i.e. smaller subunits, chosen by the respective study programmes. Each block has a corresponding workshop treating a `real-world' problem relevant to the given set of students, and these workshops form the basis of the oral exam.

Through evaluation of the comments given by students at the steering group meetings during the semesters of autumn 2021 and spring 2022, we see that students do tend to appreciate the workshop problems if they are sufficiently close to their prospective field of study. One needs to choose problems carefully, though, as problems that veer too far from the lecture contents may overwhelm students.

As we also highlight, however, performing such a restructuring is non-trivial, and there are a number of administrative challenges that one needs to keep in mind when implementing a similar structure. In particular, giving each study programme more flexibility in designing their curricula results in less flexibility when scheduling courses.

\section*{Acknowledgments}
Restructuring the calculus and linear algebra courses has involved many people from the initial conception of the idea to the full implementation. This not only includes our colleagues at the Department of Mathematical Sciences, but also teachers and researchers from departments whose students follow the course, as well as management at all levels supporting the idea.
Without their collaboration, the new structure would have been impossible.

% argument is your BibTeX string definitions and bibliography database(s)
\bibliographystyle{IEEEtran}
\bibliography{Bibliography.bib}

\vfill

\end{document}